\documentclass[12pt]{amsart}
\usepackage{amsmath}
\usepackage{amsfonts,amssymb,amsthm, txfonts, pxfonts, boxedminipage}

\DeclareMathOperator{\acosh}{arccosh}

\newtheorem{theorem}{Theorem}
\newtheorem{lemma}{Lemma}
\newtheorem{corollary}{Corollary}
\theoremstyle{remark}
\newtheorem{remark}[theorem]{Remark}

\begin{document}

\title[Continuity of volumes]{Continuity of volumes -- on a generalization of a conjecture of J.W.Milnor}
\author{Igor Rivin}
\begin{abstract}
In his paper \cite{milnorcol} J. Milnor conjectured that the volume of
$n$-dimensional hyperbolic and spherical simplices, as a function of the
dihedral angles, extends continuously to the closure of the space of
allowable angles. A prove of this conjecture was recently given by F. Luo
(\cite{luovol}\footnote{Luo's argument uses Kneser's formula \cite{kneser36}
  together some some delicate geometric estimates}). 
In this paper we give a simple proof of this conjecture, prove much sharper
regularity results, and then extend the method to apply to a large class of
convex polytopes. The simplex argument works without change in dimensions
greater than $3$ (and for spherical simplices in all dimensions), so the bulk
of this paper is concerned with the three-dimensional argument. The estimates
relating the diameter of a polyhedron to the length of the systole of the
polar polyhedron are of independent interest.\end{abstract} 
\thanks{The author would like to thank D. Laksov for his continued encouragement.}

\date{\today}

\maketitle
\section{Introduction}
Consider the set of simplices in $\mathbb{H}^n$ or $\mathbb{S}^n.$ It is
well-known that this set is parametrized by the (ordered) collection of
dihedral angles, and we may call the set of assignements of dihedral angles of
geometric simplices in $\mathbb{H}^n$ as a subset $\Omega_{\mathbb{H}^n}
\subset \mathbb{R}^{n(n+1)/2 }$ and similarly, the set of dihedral angle
assignements of of geometric simplices in $\mathbb{S}^n$ as a subset
$\Omega_{\mathbb{S}^n} \subset \mathbb{R}^{n(n+1)/2}.$  These sets are open,
since they are defined by collections of strict inequalities (which are
polynomial in the cosines of the dihedral angles).  One may then view the
volume $V$ of a simplex as a function $V$ on $\Omega_{X^n}.$ A natural
question, then, is that of boundary regularity of $V$ on $\overline{\Omega}.$ 
This question was, in fact, asked (in a weaker form of: does $V$ admit a
continuous extension to $\overline{\Omega}$) by John Milnor\footnote{Milnor
  does not attribute the conjecture to himself, so it might be
  older. W. Thurston (personal communication) opines that the conjecture is,
  in fact, due to Milnor}  in \cite{milnorcol}.  An answer to Milnor's
question was given to F. Luo in \cite{luovol}.

In this paper we begin by answering Milnor's question in the sharp form (our
method gives sharp regularity estimates, though we don't dwell too much on this). In
dimension greater than $3$ the argument goes through \emph{verbatim} to give
the same result for general convex polytopes, as long as they are of a
combinatorial class where the dihedral angles determine the volume (this
includes, eg, all \emph{simple} polytopes, that is, those with simplicial
vertex links). In dimension $3,$ things are more subtle yet, but we get the
result for various classes of convex polyhedra (those with non-obtuse dihedral
angles, by using Andreev's theorem \cite{andreev1} and for ``hyperideal''
polyhedra, by using the results of Bao-Bonahon \cite{bbpoly}). It is quite
likely analogous results should hold for hyperbolic cone manifolds. 
The omnibus form of our results on volumes is:
\begin{theorem}
Let $\Omega$ be one of the following parameter spaces:
\begin{enumerate}
\item The space of dihedral angles of polytopes of a fixed combinatorial type
  in $\mathbb{S}^n,$ such that the volume is determined by the dihedral
  angles. 
\item The space of dihedral angles of polytopes of a fixed combinatorial type
  in $\mathbb{H}^n,$ $n > 3,$ such that the volume is determined by the
  dihedral angles. 
\item The space of dihedral angles of  tetrahedra (possibly with some
  hyperideal vertices) in $\mathbb{H}^3.$ 
\item The space of dihedral angles of a fixed combinatorial type of
  \emph{hyperideal} polyhedra in $\mathbb{H}^3.$ 
\item The space of \emph{non-obtuse} dihedral angles of a fixed combinatorial
  type of polyhedra in $\mathbb{H}^3.$ 
\end{enumerate}
The the volume function on $\Omega$ extends to a $C^{0,1}$ function on $\overline{\Omega}.$
\end{theorem}

\begin{corollary} To the above, we can add: The space of dihedral angles of a
  fixed combinatorial type of \emph{ideal} polyhedra in $\mathbb{H}^3.$ 
\end{corollary} 
\begin{proof} These lie in the boundary of the set of hyperideal polyhedra, and are determined by their dihedral angles (see, eg, \cite{ann2}). In fact, the same argument works for polyhedra with some ideal and some hyperideal vertices, by the work of \cite{bbpoly}.
\end{proof}
The plan of the paper is as follows. 

In Section \ref{simplices} we give the argument for simplices. This has most
of the ingredients of the general result. 

In Section \ref{convexpoly} we recall the basic characterization results for
convex polyhedra in $\mathbb{H}^3$ in terms of their dihedral angles (or, more
generally, polar metric). 

In Section \ref{degennes} we give quantitative estimates of how close a
polyhedron is to degeneracy (boundary of the allowable dihedral angle space)
in terms of its diameter. The results (which are of independent interest) can
be loosely summarized as follows: 

\begin{theorem}
Let $P$ be a polyhedron with $N$ vertices in $\mathbb{H}^3$ of diameter $\rho
\gg 1.$ Let $M^*$ be the polar metric of $P$ (as in \cite{thes,thes0}). The
$M^*$ lies within $c_1(N) \exp(-c_2(N) \rho)$ of the boundary of the space of
admissible polar metrics, where $c_1, c_2$ are strictly positive functions of
$N$. 
\end{theorem}
The constants in the statement of the Theorem above are completely explicit,
and can be sharpened by taking into consideration finer invariants of the
combinatorics of $P$ then the number of vertices.

\section{A simple proof for simplices (among other things)}
\label{simplices}

In dimension $2,$ the result follows immediately from Gauss' formula, which
states that area is a linear function of the angles, so we will only discuss
dimensions $3$ or above. 

The simple proof relies on the \emph{Schl\"afli differential equality} (see 
\cite{milnorcol}, which states that in a space of constant curvature $K,$ and
dimension $n$ the volume of a smooth family of polyhedra $P$ satisfies the
differential equation: 

\begin{equation}
\label{schlafli1}
K d V(P) = \dfrac{1}{n-1} \sum_F V_{n-2}(F) d \theta_F,
\end{equation}

where the sum is over all codimension-$2$ faces, $V_{n-2}$ is the $n-2$ dimensional volume of $F,$ and $\theta_F$ is the dihedral angle at $F.$

Another way of writing the Schl\"afli formula is:

\begin{equation}
\label{schlafli2}
K \dfrac{\partial V(P)}{\partial \theta_F} = V_{n-2}(F).
\end{equation}

This is the form we will use. 

The first observation is that $V_{n-2}(F)$ is bounded by a constant
(dimensional for  $\mathbb{S}^n,$ , depending on the number of vertices of
$F,$ but nothing else, and  if $n \geq 4,$  in $\mathbb{H}^n.$).  

This immediately shows the continuity of volume for all $\mathbb{S}^n,$ and
for $\mathbb{H}^n,$ whenever $n \geq 4.$ 

We are left with dimension $3.$ All we really need is noting that the partial
derivatives of $V$ with respect to the dihedral angles develop at worst
logarithmic singularities as we approach the frontier of
$\overline{\Omega_{\mathbb{H}^3}}$ -- this result suffices by the following
form of the Sobolev Embedding Theorem (this is \cite[Theorem 7.26]{gilbrtru}): 
\begin{theorem}
\label{sobemb}
Let $\Omega$ be a $C^{0,1}$ domain in $\mathbb{R}^n.$ Then,
\begin{itemize}
\item{(i)} If $k p < n,$ the space $W^{k,p}(\Omega)$ is continuously imbedded in $L^{p^*}(\Omega),$ where $p^* = np/(n - kp),$ and compactly imbedded in $L^q(\Omega)$ for any $q < p^*.$
\item{(ii)} If $0 \leq m < k - \frac{n}{p} < m+1,$ the space $W^{k, p}$ is continously embedded in $C^{m, \alpha}(\overline{\Omega}), $ $\alpha=k -n/p - m,$ and compactly embedded in $C^{m, \beta}(\overline{\Omega})$ for any $\beta < \alpha.$
\end{itemize}
\end{theorem}
Here, the Sobolev space $W^{k, p}$ is the space of functions whose first $k$ (distributional) derivatives are in $L^p.$

In our case, we know that the domain $\Omega$ is bounded, convex ``curvilinear polyhedral'' (hence $C^{0, 1}$) domain,  volume is a bounded function, and we assume that the gradient grows logarithmically as we approach the boundary. This implies that $V$ is in $W^{1, p}$ for \emph{all} $p > 0,$ so we get the following corollary:
\begin{corollary}
Volume is in $C^{0, \alpha}(\overline{\Omega})$ for any $\alpha < 1.$
\end{corollary}

The logarithmic growth of diameter of the simplex as a function of the distance to $\partial \Omega$ can be shown in a completely elementary way using Eq. \eqref{schlafli2} and elementary reasoning about Gram matrices, as follows:

Let $G$ be ``angle Gram matrix'' of a simplex $\Delta$, that is, $G_{ij} = - \cos\theta_{ij},$ where $\theta_{ij}$ is the angle between the $i$-th and the $j$-th face.  Let $S$ be the matrix whose columns are the normals to the faces of $\Delta$ (all the computations take place in Minkowski space, and we use the hyperboloid model of $\mathbb{H}^n.$ It is immediate that $G = S^t S.$

Let now $W$ be the matrix whose columns are the (possibly scaled) vertices of $\Delta.$ $W$ satisfies the equation $S^t W = I,$ and to get the vertices to lie on the hyperboloid $\langle x, x\rangle = -1$ we must rescale in such a way that the squared norms of the columns of $W$ become $-1.$ Call the scaled matrix $W_s.$ Since the usual ``length'' Gram matrix $G^*$ of $\Delta$S can be written as $W_s^tW_s,$ and $G^*_{ij} = - \cosh(d(v_i, v_j)),$ a simple computation using Cramer's rule gives:
\[
\cosh d(v_i, v_j) = \dfrac{c_{ij}}{\sqrt{c_{ii} c_{jj}}},
\]
where $c_{ij}$ is the $ij$-th cofactor of $G.$ (see \cite{muraushi} for many related results).

It  follows that the distances between the vertices (which are the lengths of the edges, which are the faces of codimension $2.$) behave as $|\log\ c_{ii}|.$ Since the cofactors are polynomial in the cosines of the angles, we are done.

It should be noted that this argument works \emph{mutatis mutandis} for \emph{hyperideal} simplices, or simplices with some finite and some hyperinfinite vertices..

\section{Convex polytopes}
\label{convexpoly}
For arbitrary convex polytopes in dimension $n > 3,$ (and convex \emph{spherical} polytopes in all dimensions)  the proof given in Section \ref{simplices} goes through without change, with the one proviso that it is not currently known whether the volume of a polytope is determined up to congruence by its dihedral angles. Such a uniqueness result \emph{is} conjectured (indeed, it is conjectured that a polytope is determined up to congruence by the dihedral angles), and is easy to prove for \emph{simple} polytopes -- those with simplicial links of vertices -- this follows in arbitrary dimension from the corresponding result in $3$ dimensions. For the $3$-dimensional result see, eg, \cite{thes0,thes}.  In the rest of this paper we will be concerned with the most interesting case of convex polyhedra in $\mathbb{H}^3.$

First, we recall the following result of \cite{thes,thes0}:
\begin{theorem}[\cite{thes,thes0}]
\label{annelk}
A metric space $(M, g)$ homeomorphic to $\mathbb{S}^2$ can arise as the Gaussian image $G(P)$ of a compact convex polyhedron $P$ in $\mathbb{H}^3$ if and only if the following conditions hold:
\begin{itemize}
\item{(a)} The metric $g$ has constant curvature $1$ away from a finite collection of cone points.
\item{(b)} The cone angle at each $c_i$ is greater than $2\pi.$
\item{(c)} The lengths of closed geodesics of $(M, g)$ are all strictly greater than $2\pi.$
\end{itemize}
\end{theorem}
Assume now, for simplicity, that the polyhedron $P$ is simple (as pointed out above, this also allows us to assume that the geometry is determined by the dihedral angles). In that case, the space of admissible metrics $\Omega_P$ (as per Theorem \ref{annelk}) is parametrized by the exterior dihedral angles (the cell decomposition dual to that of $P$ gives a triangulation of the Gaussian image, and the (exterior) dihedral angles are the lengths of edges of the triangulation.) Theorems \ref{degenthm},\ref{quasig} immediately imply the following:

\begin{theorem}
\label{loggrowth}
There exists a constant $L_0,$ such that 
the maximal length $\ell_P$ of an edge of $P$ is bounded as follows:
\[
\ell_P \leq \max(L_0,  - 2N \log(d(P, \partial \overline{\Omega}_P)/12N)),
\]
where $N$ is the number of vertices of $P.$
\end{theorem}
\begin{proof}
Assume the contrary. Then, there exists a sequence of polyhedra $P_1, \dots, P_n, \dots$ with diameter $\rho(P_i) \geq \ell){P_i}$ going to infnity, which are farther than $12 N  \exp(-\rho/2N).$ By choosing a subsequence, we may assume that there is a \emph{fixed} cycle of faces $F_1, \dots, F_k$ of $P,$ such that the sum of dihedral angles along the edges $e_i = F_i \cap F_{i+1}$ is smaller than $2\pi + 12 N  \exp(-\rho/2N),$ (by Theorem \ref{degenthm}) and which are a $4N\exp(-2\rho)$ quasigeodesic (by Theorem \ref{quasig}). Since the limit point of the $P_i$ is not in $\Omega_P$ (by Theorem \ref{annelk}),
the result follow.
\end{proof}
The following corollary is immediate (by Schl\"afli, see Section \ref{simplices}):
\begin{corollary}
\label{lpest}
The volume is in $W^{1, p}(\Omega_p)$ for \emph{all} $p > 0.$
\end{corollary}
We now have almost enough to show that volume extends to $\overline{\Omega}_p,$ except for the slight matter of not having the required (by Theorem \ref{sobemb}) regularity result for $\partial{\Omega}_P.$ Such a result seems quite non-trivial, since the length of the shortest closed geodesic is a rather badly behaved quantity. However, there are two very important special cases where  $\Omega_P$ is actually a convex polytope, so this problem is finessed. The first is when we restrict the dihedral angles to be non-obtuse. This case is covered by 
\begin{theorem}[Andreev's theorem (\cite{andreev1})]
Let $C$ be an abstract polyhedron with more than $4$ faces, and suppose that \emph{non-obtuse} angles $\alpha_i$ are given corresponding to each edge $e_i$ of $C.$ There is a compact convex hyperbolic polyhedron whose faces realize $C$ with dihedral angle $\alpha_i$ at each edge $e_i$ if and only if the following conditions hold:
\begin{enumerate}
\item For each edge $e_i,$ $0 < \alpha_i \leq \pi/2.$
\item For any three distinct edges $e_1, e_2, e_3$  meeting at a vertex, $\alpha_1 + \alpha_2 + \alpha_3 > \pi.$
\item Whenever $\Gamma$ is a prismatic $3$-circuit intersecting edges $e_1, e_2, e_3,$ then $\alpha_1 + \alpha_2 + \alpha_3 < \pi.$
\item Whenever $\Gamma$ is a prismatic $4$-circuit intersecting edges $e_1, e_2, e_3, e_4,$ then
$\alpha_1 + \alpha_2 + \alpha_3 + \alpha_4 <2\pi.$
\item Whenever there is a four-sided face bounded by edges $e_1, e_2, e_3, e_4$ (enumerated successively), with edges $e_{12}, e_{23}, e_{34}, e_{41}$ entering the four vertices, then:
\begin{gather}
\alpha_1 + \alpha_3 + \alpha_{12} + \alpha_{23} + \alpha_{34} + \alpha_{41} < 3\pi,\\
\alpha_2 + \alpha_4 + \alpha_{12} + \alpha_{23} + \alpha_{34} + \alpha_{41} < 3\pi.
\end{gather}
\end{enumerate}
\end{theorem}

The other special case is that of hyperideal polyhedra:
\begin{theorem}[Bao-Bonahon, \cite{bbpoly}]
Let $\sigma$ be a cell decomposition of $S^2,$ and let $w: \sigma_1 \rightarrow (0, \pi)$ be a map on the set of edges of $\sigma.$ Then there exists a hyperideal polyhedron with combinatorics given by $\sigma$ and exterior dihedral angles given by $w$ if and only if:
\begin{enumerate}
\item
The sum of the values of $w$ on each circuit in $\sigma_1$ is not smaller than $2\pi,$ and is strictly greater if the circuit is non-elementary. 
\item The sum of values of $w$ on each simple path in $\sigma_1$ is strictly larger than $\pi.$
\end{enumerate}
The hyperideal polyhedron is then unique.
\end{theorem}
\section{Degeneration estimates}
\label{degennes}
The results of this section are a quantitative version of the results of the compactness results of \cite{thes,thes0}. 
First, some key lemmas. The general setup will be as follows: $L$ is a geodesic in $\mathbb{H}^3,$  $t$ is a real number (generally large) and $P, P^-, P^+$ are three planes, all orthogonal to $L,$ and such that $d(P, P^-) = d(P, P^+) = t,$ and $d(P^-, P^+) = 2 t.$ We denote $x_0 = L \cap P.$

In the sequel, we use the hyperboloid model of $\mathbb{H}^3,$ where $\mathbb{H}^3$ is represented by the set $\langle x, x \rangle = -1;$ $x_0 > 0,$ in the $\mathbb{R}^4$ equipped with the scalar product $\langle x, y \rangle = - x_1 y_1 + \sum_{i=2}^4 x_i y_i.$ The reader is referred to \cite{wptbook} (as well as \cite{thes}) for the details (which will be used below).

Returning back to our setup, we can assume, without loss of generality, that
\[
x_0=\begin{pmatrix}1\\0\\0\\0\end{pmatrix},
\]
that 
\[
P^\perp = \begin{pmatrix}0\\1\\0\\0\end{pmatrix},
\]
and hence, that $P^+ = \phi(t) P,$ while $P^- = \phi(-t) P,$ where
\[
\phi(r) = \begin{pmatrix} \cosh(r) & \sinh(r) & 0 & 0\\
                                           \sinh(r) & \cosh(r) & 0 & 0\\
                                           0 & 0 &1 &0\\
                                           0 & 0 & 0 &1
                                           \end{pmatrix}.
                                           \]
                                           
Since $\phi(r)$ is symmetric, it follows that \[P^{+\perp} = \phi(t) P^\perp = \begin{pmatrix} \cosh(t) \\ \sinh(t) \\ 0 \\0\end{pmatrix},\] while 
\[P^{-\perp} = \phi(t) P^\perp = \begin{pmatrix} \cosh(t) \\ -\sinh(t) \\ 0 \\0\end{pmatrix},\]
\begin{lemma}
\label{anglemma}
 Let $Q$ be a plane in $\mathbb{H}^3$ which intersects both $P^-$ and $P^+.$ Then, there exists $t_0,$ such that $Q$ intersects $P,$ and the cosine of the  angle $\alpha$ of intersection satisfies $|\cos(\alpha)| <3 e^{-t},$ as long as $t > t_0.$ The number $t_0$ can be picked \emph{independently} of $Q.$
\end{lemma}
 
\begin{proof}
Let the unit normal $Q^\perp$ to $Q$ be $Q^\perp=\begin{pmatrix}a\\b\\c\\d\end{pmatrix}.$  
Since two planes intersect if and only if the scalar product of their unit normals is less than $1$ in absolute value, we have, from the hypotheses of the lemma and the description of the unit normals to $P^-$ and $P^+$ above that:
\begin{gather}
|a \cosh(t) + b \sinh(t)| < 1 \\
|a \cosh(t) - b \sinh(t)| < 1.
\end{gather}
Squaring the two inequalities, and adding them together we obtain:
\[
a^2 \cosh^2(t) + b^2 \sinh^2(t) < 1.
\]
Since, under the hypotheses of the lemma, $\min(\cosh(t), \sinh(t)) > e^t/3,$ it follows that
\[
a^2 + b^2 < 3/e^t,
\]
and so $\max(a, b) < 3e^{-t}.$
Now, the cosine of the angle between $Q$ and $P$ equals $\langle Q^\perp, P^\perp\rangle = b,$ so the result follows.
\end{proof}

\begin{remark}
The constant $3$ is far from sharp (especially for larger $t$).
\end{remark}

\begin{lemma} 
\label{distlem}
There exists a $t_0,$ such that if  $M$ is a line in $\mathbb{H}^3$ which intersects both $P^-$ and $P^+,$ then $M$ intersects $P,$ and $\cosh(d(P\cap M, x_0)) < 4 e^{-2t}+1, $ as long as $t > t_0.$
\end{lemma}

\begin{proof}
Assume that 
$M\cap P^+ = \phi(t)  p_1,$ and $M\cap P^- = \phi(-t) p_1,$ where $p_{1, 2} \in P.$ (This is always possible, since $P^+ = \phi(t) P,$ $P^- = \phi(-t).$) The intersection of $M$ with $P$ is then given by 
\[
M\cap P = \dfrac{x (M \cap P^+) + y(M \cap P^-)}{\|x (M \cap P^+) + y(M \cap P^-)\|},
\]
where $x$ and $y$ are chosen so that the linear combination is actually in $P,$ or, in other words, the second coordinate of the linear combination vanishes. We abuse notation above by writing
$\|Z\| = \sqrt{-\langle Z, Z\rangle}.$

Let us now compute. Set (for $i=1, 2$)
\[
p_i = \begin{pmatrix} a_i\\0\\c_i\\d_i\end{pmatrix}.
\]
It follows that 
\[
M\cap P^+ = \begin{pmatrix}a_1 \cosh(t) \\ a_1 \sinh(t)\\ c_1 \\ d_1\end{pmatrix}.
\]
while
\[
M \cap P^- = \begin{pmatrix}a_2 \cosh(t) \\ - a_2 \sinh(t) \\ c_2 \\d_2 \end{pmatrix}.
\]
It follows that we can choose $x = 1/(2a_1),$ $y = 1/(2a_2),$
so that
\[
m = x M\cap P^+ + y M\cap P^- = \begin{pmatrix} \cosh(t) \\0 \\ \frac{1}{2}( c_1/a_1 + c_2/a_2) \\ \frac{1}{2} ( d_1/a_1 + d_2/a_2) \end{pmatrix}.
\]
It follows that 

\begin{multline}
- \cosh(d(M \cap P, x_0)) = \left\langle \dfrac{m}{\|m\|}, x_0\right\rangle = \\ -\dfrac{\cosh(t)}{\sqrt{\cosh^2(t) - 1/4\left((c_1/a_1 + c_2/a_2)^2 + (d_1/a_1 + d_2/a_2)^2\right)}}.
\end{multline}

Since $c_i^2 + d_i^2 + 1 = a_i^2,$ for $i=1,2$ it follows that $|c_i/a_i| < 1,$ and similarly $|d_i/a_i| < 1,$ so that 
\[
\cosh^2(t) \geq \cosh^2(t) - 
 1/4\left((c_1/a_1 + c_2/a_2)^2 + (d_1/a_1 + d_2/a_2)^2\right) > \cosh^2(t) -2.
 \]
 It follows that 
 \[
\cosh(d(M \cap P, x_0))\leq \dfrac{1}{\sqrt{1-2/\cosh^2(t)}},
 \]
 and the assertion of the lemma follows by elementary calculus.
\end{proof}

\begin{lemma}
\label{spherical}
Let $T$ be a spherical triangle with sides $A,B,C$ and (opposite) angles $\alpha, \beta, \gamma.$ Suppose that $|\cos(\beta)| < \epsilon \ll 1,$ $| \cos(\gamma)| < \epsilon \ll 1.$ Then 
$|\alpha - A| < 2 \epsilon|.$
\end{lemma}

\begin{proof} The spherical Law of Cosines states that:
\[
\cos(A) = \dfrac{\cos(\alpha)+ \cos(\beta) \cos(\gamma)}{\sin(\beta) \sin(\gamma)}.
\]
It follows that 
\[
\cos(A) - 2\epsilon^2 \leq \cos(A)(1-\epsilon^2) - \epsilon^2  \leq \cos(\alpha) \leq \cos(A) + \epsilon^2.
\]
The assertion of the lemma follows immediately.
\end{proof}

\begin{corollary}
\label{spherecor}
Let $F_1$ and $F_2$ be two planes intersecting at a dihedral angle $\alpha,$ with both $F_1$ and $F_2$ intersecting a third plane $P,$ at angles whose cosines are smaller than $\epsilon.$ Let $A$ be the angle between $F_1 \cap P$ and $F_2 \cap P.$ Then $|\alpha - A|< 2\epsilon.$
\end{corollary}

\begin{proof}
Apply Lemma \ref{spherical} to the link of the point $F_1 \cap F_2 \cap P.$
\end{proof}

\begin{lemma}
\label{circleest}
Let $V$ be a convex polygon in the hyperbolic plane $\mathbb{H}^2,$ such that all the vertices of $V$ lie within a distance $r$ of a certain point $O.$ Then, the sum of the exterior angles of $V$ is smaller than $2\pi \cosh(r).$
\end{lemma}

\begin{proof}
The area of a disk of radius $r$ in $\mathbb{H}^2$ equals $4\pi \sinh^2(r/2) = 2\pi (\cosh(r) -1)$ (see \cite{vinberg}). 
Since $V$ is contained in such a disk, its area is at most $2\pi(\cosh(r) -1),$ and since the area of $V$ equals the difference between the sum of the exterior angles and $2\pi,$ the statement of the lemma follows.
\end{proof}

Now we are ready to show the following:
\begin{theorem}
\label{degenthm}
Let $X$ be a convex polyhedron with $N$ vertices  in $\mathbb{H}^3$ of diameter $\rho \gg 1.$ Then, there exists a cyclic sequence of faces $F_1, \dotsc, F_k = F_1,$ with $F_i$ sharing an edge $e_i$  with $F_{i+1}$ (indices taken $\mod k$) so that the sum of exterior dihedral angles at $e_1, \dotsc, e_k$ is smaller than $2\pi + 12 N  \exp(-\rho/2N).$
\end{theorem}

\begin{proof}
Take a diameter $D$ of $X$  of length $\rho,$ place points $p_1, \dotsc, p_N$ equally spaced on $D.$ By the pigeonhole principle, one of the segments $p_ip_{i+1}$ contains no vertices of $X.$ Let $x_0$ be the midpoint of the segment $p_ip_{i+1}.$ Construct planes orthogonal to $D$ at $x_0$ ($P$) and $p_i$ ($P^-$), and at $p_{i+1}$ ($P^+$). Let $t=\rho/(2N).$ The portion of $X$ contained between $P^-$ and $P^+$ is a polyhedral cylinder, consisting of faces $F_1, \dots, F_k.$ By Lemma \ref{distlem}, the intersection of $X$ with $P$ is a polygon $\mathcal{P},$ whose sum of exterior angles is at most $2\pi(4\exp(-2t) + 1),$ and so by Corollary \ref{spherecor}, combined with Lemma \ref{anglemma}, the sum of the dihedral angles corresponding to pairs $F_i F_{i+1}$ is at most $2\pi(4 \exp(-2t) + 1) + 6 k \exp(-t).$ Since $k$ is no greater than the number of faces of $X,$ which, in turn, is at most $2N -4.$
\end{proof}

\begin{theorem}
\label{quasig}
With notation as in Theorem \ref{degenthm}, the faces $F_1, \dotsc, F_k$ form a  curve in the Gaussian image of $X$ with geodesic curvature not exceeding $3k\exp(-\rho/N)).$
\end{theorem}
\begin{remark} The reader is referred to \cite{thes0,thes} for a more thorough discussion of geodesics on spherical cone manifold, but suffice it to say that the contribution of the face $F_i$ to the geodesic curvature is $0$ if the two edges are (hyper)parallel, and equal to the angle of intersection if they intersect.
\end{remark}
\begin{proof}
Let $e_1$ and $e_2$ be the two edges of $F.$ If $e_1$ and $e_2$ do not intersect, there is nothing to prove (by the remark above. If they do intersect at a point $C,$ note that $C$ is at a distance at least $\rho/2N$ from $x_0,$ while the intersections $A$ and $B$ of $e_1$ and $e_2$ with $P$ are at most 
$\acosh(4\exp(\-rho/N)+1) \approx \sqrt{8} \exp(-\rho/2N)$ away from $x_0,$ and so at most (for large $\rho$) $6\exp(-\rho/2N)$ away from each other. We will only use the (much cruder) estimate $\cosh(AB) \leq 2.$ Now, apply the hyperbolic law of cosines to the triangle $ABC,$ to get:

\begin{multline}
1\geq \cos(\gamma) = \dfrac{-\cosh(AB) + \cosh(AC)\cosh(BC)}{\sinh(AC) \sinh(BC)} \geq \\
1- \dfrac{\cosh(AB)}{\sinh(AC)\sinh(BC)} \geq 1- \dfrac{2}{\sinh(AC)\sinh(BC)} \geq 1 - 8 \exp(-\rho/2N).
\end{multline}
The estimate now follows.

\end{proof} 

\begin{remark}
The argument above is easily modified to show that the curve dual to $F_1, \dotsc, F_k$ has small geodesic curvature viewed as a curve in $\mathbb{S}_1^3,$ and not just in $X^*.$
\end{remark}

\bibliographystyle{plain}
\bibliography{curves,rivin,opt,matrix}
\end{document}